\documentclass[a4paper,12pt]{article}
\usepackage[colorlinks=true,urlcolor=blue]{hyperref}
\usepackage{amsmath}
\usepackage{amssymb}
\usepackage{amsthm}
\usepackage{fancyhdr}
\begin{document}
\parindent 0cm

\textbf{\large{About the Algebraic Solutions of Smallest Enclosing\\Cylinders Problems}}
\vskip 1cm

\textbf{Michel Petitjean}
\vskip 1cm

\small{
MTi, INSERM UMR-S 973, University Paris 7

35 rue H\'el\`ene Brion, 75205 Paris Cedex 13, France.

E-mail: petitjean.chiral@gmail.com

http://petitjeanmichel.free.fr/itoweb.petitjean.html
}

\vskip 1cm
\textbf{Abstract}.
Given $n$ points in Euclidean space $E^d$, we propose an algebraic algorithm
to compute the best fitting $(d-1)$-cylinder.
This algorithm computes the unknown direction of the axis of the cylinder.
The location of the axis and the radius of the cylinder are deduced analytically
from this direction.
Special attention is paid to the case $d=3$ when $n=4$ and $n=5$.
For the former, the minimal radius enclosing cylinder is computed algebrically
from constrained minimization of a quartic form of the unknown direction of the axis.
For the latter, an analytical condition of existence of the circumscribed cylinder
is given, and the algorithm reduces to find the zeroes of
an one unknown polynomial of degree at most $6$.
In both cases, the other parameters of the cylinder are deduced analytically.
The minimal radius enclosing cylinder is computed analytically
for the regular tetrahedron and for a trigonal bipyramids family with a symmetry axis of order 3.



\parindent 0cm
\vskip 1cm
\textbf{1. Introduction}
\vskip 0.5cm
\parindent 0.5cm

Let $E^d$ be the $d$-dimensional Euclidean space,
$U^{d-j}$ be a $(d-j)$-dimensional affine subspace,
and $\rho$ a non negative real number.
The set of points $x$ of $E^d$ lying at distance smaller or equal to $\rho$
from their orthogonal projection on $U^{d-j}$ is called a $j$-cylinder.
We consider $n$ data points in $E^d$.
A $j$-cylinder $C$ containing $n$ points is called a circumscribed $j$-cylinder
if all $n$ points are lying on the boundary of $C$.
We consider further only $(d-1)$-cylinders in this paper.
A $(d-1)$-cylinder is defined with an axis, i.e. a unit vector $u$
defining the direction of the axis and a point $c$ of $E^d$ locating the axis,
and a radius $\rho$.

Given $n$ data points $x_i, i=1,...,n$ in $E^d$,
several cylinder problems are defined.

\textit{Best fitting cylinder problem}:
Find the $(d-1)$-cylinder $(u,c,\rho)$ minimizing
the standard deviation of the population of the $n$ squared distances
between the $x_i$ and their projections on the axis,
$\rho^2$ being the mean of this population.
If it happens that the minimal standard deviation is null,
all distances are equal to $\rho$ and
the resulting cylinder is circumscribed to the $n$ points.

\textit{Smallest enclosing cylinder problem}:
Find the $(d-1)$-cylinder $(u,c,\rho)$ of minimal radius $\rho$
enclosing the $n$ points.

\textit{Circumscribed cylinder problem}:
Find the $(d-1)$-cylinder $(u,c,\rho)$ circumscribed to the $n$ points,
and if it is not unique, exhibit the one of minimal radius.
The circumscribed cylinder may not exist.

The circumscribed cylinder problem in $E^3$ seems to have appeared in the literature
in 1977 \cite{Bottema1977}, but practical algorithms appeared much later \cite{Devillers2003},
at the occasion of the analysis of circular cylinders through four or five points.
Unfortunately, the practical computation of solutions is not trivial and
requires the use of specialized packages. Recently, a numerical solver based
on second order cone programming methods has been devised to compute
smallest enclosing cylinders \cite{Watson2006}. Bounds for the radius of
minimal enclosing $d$-simplices have been given \cite{Brandenberg2006}, and
the complexity of the smallest enclosing cylinder problem has been investigated 
for metrology applications \cite{Schomer2000}.

In this paper we propose simple algebraic algorithms to solve the three cylinders problems,
which can be implemented without the use of sophisiticated packages.
The motivation comes from recent molecular modeling studies
in which the shape of many small molecules is expected to be
suitably fitted with cylinders \cite{Meslamani2009}.
In this context, $d=3$, and the smallest enclosing cylinder problem
is the basic one we need to solve. We face to small $n$ values,
so that we are interested by the algebraic aspects rather than by the computational complexity.
Solving the best fitting cylinder problem is useful to introduce our approach,
and solving the circumscribed cylinder problem is part of our full algorithm.

\parindent 0cm
\vskip 1cm
\textbf{2. The best fitting cylinder problem}
\vskip 0.5cm
\parindent 0.5cm

The usual scalar product of $x$ and $y$ is $x'y$,
the quote denoting a vector of matrix transposition.
The norm of $x$ is $\Vert x\Vert=(x'x)^{1/2}$.
We define $X$ as the rectangular array of $n$ lines and $d$ columns
containing $x_i'$ at the line $i$, $i=1,...,n$.
The squared distance between $x_i$ and its projection on the axis is
$\Delta_i = (x_i-c)'(x_i-c) - ((x_i-c)'u)^2$.
The mean of the population of the $\Delta_i$ is $\bar{\Delta}$,
and the variance to be minimized is
$\frac{1}{n}\sum(\Delta_i-\bar{\Delta})^2$.

Any point $c$ lying on the axis of the cylinder can be used.
Thus we decide to retain $c$ as the projection of the origin on the axis, i.e. $c'u=0$.
We are left with the minimisation of a polynomial of $u$ and $c$,
subject to the orthogonality condition $c'u=0$ and to the normalization $u'u=1$.
Using the Lagrange multipliers K and L, the function to be minimized is

\parindent 2cm
\vskip 0.5cm
                       $Min_{\lbrace c,u\rbrace}F_{cu}$ \hfill (1)

                       $F_{cu}=\frac{1}{n}\sum(\Delta_i-\bar{\Delta})^2 - Kc'u - L(u'u-1)$
\vskip 0.5cm
\parindent 0.5cm

For clarity we assume further than the mean of the $n$ points is translated to the origin.
Setting $V_i=x_i x_i'$, the inertia matrix is $T=\sum V_i$ and the covariance matrix is $V=T/n$.
We also further assume that the set of the $n$ points is not subdimensional,
i.e. $X$ is of full rank $d$, and $T=X'X$ and $V$ are invertible.
The identity matrix of rank $d$ is $I$,
we define $B_i=I\cdot Tr(V_i-V) - (V_i-V)$.
Then, $\Delta_i=Tr(Vi)-2c'x_i+c'c-u'V_iu$, and the expression of $F_{cu}$ expands as

\parindent 2cm
\vskip 0.5cm
               $F_{cu}=\frac{1}{n}\sum(u' B_i u - 2c'x_i)^2 - Kc'u - L(u'u-1)$ \hfill (2)
\vskip 0.5cm
\parindent 0.5cm

The components of $u$ and $c$ are the $2d$ unknowns in (2). The unknown radius
is $\rho^2=\bar{\Delta}$. It is computed from $c$ and $u$.
The gradient of $F_{cu}$ relative to $c$ is

\parindent 2cm
\vskip 0.5cm
               $G_c=-\frac{4}{n}\sum(x_i u' B_i u) +8Vc -Ku$ \hfill (3)
\vskip 0.5cm
\parindent 0.5cm

The stationary point is such that $G_c=0$, so that

\parindent 2cm
\vskip 0.5cm
               $8c=V^{-1} \lbrack Ku +\frac{4}{n}\sum(x_i u' B_i u) \rbrack$ \hfill (3)
\vskip 0.5cm
\parindent 0.5cm

The Lagrangian $K$ is expressed from setting $u'c=0$ in (3).
Then, the analytical expression of the optimal $c$ value can be used in (2), so that the minimization
in (1) reduces to a non convex optimization problem in which there are only $d$ unknowns, i.e.
the components of $u$. Standard algebraic solvers can be used \cite{Gill2008}.
Moreover, starting points are random unit vectors $u$, which in fact are $d-1$ independant parameters values.
This is beneficial in term of global optimization cost.
E.g., for $d=3$, a reasonable amount of random initial unit vectors should suffice
to locate the global minimum of $F_{cu}$, and there is no need of a huge of initial $6$-tuples.

The case where the solution is $F^*_{cu}=0$ is of special interest.
In this situation, all the $n$ points are lying at distance $\rho$ from the axis,
and we face to a circumscribed cylinder problem.
We have $n$ equations $\Delta_i=\bar{\Delta}$, i.e.

\parindent 2cm
\vskip 0.5cm
               $u' B_i u = 2c'x_i$ \hfill (4)
\vskip 0.5cm
\parindent 0.5cm

Only $n-1$ are independant equations because $\sum B_i=0$ and $\sum x_i=0$.
Adding to them the two equations $c'u=0$ and $u'u=1$, we get a system of $n+1$ independant equations
of $2d$ unknowns $c$ and $u$. This system is determined when $n=2d-1$.
When $d=3$, we retrieve that at most a finite number of cylinders is expected to pass through $5$ points
in general position. More will be said further.

\parindent 0cm
\vskip 1cm
\textbf{3. The smallest enclosing cylinder problem}
\vskip 0.5cm
\parindent 0.5cm

A cylinder is a convex set, so that the smallest cylinder enclosing $n$ points is
the one which encloses the vertices of the convex hull of the $n$ points.
Thus, deleting the points interior to the hull before computing the smallest enclosing cylinder
can drastically reduce the computational cost.
Convex hull algorithms in $2$, $3$, and more dimensions have been described
\cite{Edelsbrunner1987,Preparata1985}.

When $d=3$, the smallest enclosing cylinder problem is here the minimal radii enclosing cylinder problem,
which should not be confused with the minimal height enclosing cylinder problem.
This latter is a special case of the smallest enclosing $1$-cylinder problem
although the former is a smallest enclosing $(d-1)$-cylinder problem.
The minimal height enclosing cylinder is defined by the closest enclosing parallel planes,
and we know from proposition 3 in \cite{Brandenberg2006} that these latter planes contain
together at least one 4-tuple of vertices of the convex hull.
The radius of the minimal height enclosing cylinder is then the radius of the smallest circle
containing the projections of the points on one of the closest enclosing parallel planes.

We return back to the minimal radii enclosing cylinder problem.
Owing to the result given at the end of section 2,
the smallest enclosing cylinder must be seeked in the set of
minimal radius cylinders circumscribed to $k$ points,
$k$ taking successively increasing values from $1$ to $2d-1$.
The smallest $k$ value for which there is at least one minimal radius circumscribed cylinder
enclosing all $n$ points is retained. If there are several such circumscribed
cylinders, the smallest radius one is retained. It may be not unique.

When $d=3$, we look for minimal radius cylinders circumscribed to
successively $1$, $2$, $3$, $4$ and $5$ points.
Cases $k=1$ and $k=2$ are degenerate.
The case $k=3$ is in fact subdimensional: find the two closest lines
enclosing $3$ points. Still from proposition 3 in \cite{Brandenberg2006}
we deduce that the associated minimal radius is half of the smallest of the
three heights of the triangle defined by the three points.
Only the cases $k=4$ and $k=5$ are non trivial.
We consider them hereafter, assuming $d=3$.

\parindent 0cm
\vskip 1cm
\textbf{4. The smallest cylinder passing through 4 points}
\vskip 0.5cm
\parindent 0.5cm

Here we set $n=4$.
For convenience, we set $\gamma=2c$ and we define the vector $b$ having four components
$b_i=u' B_i u$. The four equations in (4) are rewritten in matricial form $X\gamma=b$,
and because we already assumed $X$ to be of full rank, we express $\gamma$ from $u$ in (5).

\parindent 2cm
\vskip 0.5cm
               $\gamma=T^{-1}X'b$ \hfill (5)
\vskip 0.5cm
\parindent 0.5cm

The squared radius $\rho^2=\bar{\Delta}=Tr(V)-u'Vu+c'c$
is to be minimized under the constraints $u'c=0$ and $u'u=1$.
Setting $W=I\cdot Tr(T)-T$ and using the Lagrangians $K$ and $L$,
the objective function of $u$ to be minimized is

\parindent 2cm
\vskip 0.5cm
               $F=u'u \cdot u'Wu+\gamma'\gamma-2K\gamma'u+L(1-u'u)$ \hfill (6)
\vskip 0.5cm
\parindent 0.5cm

The case of the regular tetrahedron is solved analytically.
\vskip 0.5cm

\textbf{Theorem 1.}
\textit{(a) There are three absolute minimal radius circumscribed cylinders to the regular tetrahedron.
The ratio of the minimal radius to the length of the edge is $1/2$.
The axis of the three cylinders intersect at the center of the tetrahedron
and are mutually orthogonal, following the directions parallel to the edges of the cube
having for vertices the four ones of the tetrahedron plus the four ones
got by reflection of these latter through the center.
(b) There are six absolute maximal radius circumscribed cylinders to the regular tetrahedron.
The ratio of the maximal radius to the length of the edge is $3\sqrt{2}/8$.
The axis of the six cylinders following the directions parallel to the edges of the tetrahedron}.
\vskip 0.5cm

\textit{Proof.} We use:
$x_1'=(1,1,1)$, $x_2'=(1,-1,-1)$, $x_3'=(1,-1,1)$, $x_4'=(-1,-1,1)$.

Part (a): the center is at the origin and the covariance matrix $V$ is the identity matrix $I$.
The squared radius to be minimized is $\rho^2=Tr(V)-u'Vu+c'c=2+c'c$,
under the constraints $u'c=0$ and $u'u=1$, with $\gamma=2c=T^{-1}X'b$.
The absolute minimum at $c=0$ and $\rho^2=2$ will be valid
under the condition that the four equations $u'B_iu=\gamma'x_i=0$ are indeed satisfied
for some unit vector $u$. Then, $u'B_iu=Tr(V_i-V)-u'(V_i-V)u$,
i.e. $u'B_iu=x_i'x_i-3-(u'x_i)^2+1=1-(u'x_i)^2$.
The solutions $u$ are such that $u'x_i=1$ or $u'x_i=-1$, which is written
in matricial form $Xu=\delta$, where $\delta$ is a vector with each of its four component
taking independantly the value $1$ or the value $-1$, and $\delta'\delta=4$.
So, $u=T^{-1}X'\delta=X'\delta/4$ and $u'u=1$.
We deduce that the set of solution unit vectors $u$
reduces to the set of the canonical base vectors and their opposite.
The rest of the proof of part (a) is trivial.

Part (b): the squared radius to be maximized is still $\rho^2=2+c'c$.
Noticing that $XX'/4$ is the centering operator $I-11'/n$,
where $I$ is now the identity matrix of size $n=4$ and
$\mathbf{1}$ is the vector having its four components equal to 1,
$c'c=b'XT^{-2}X'b/4=b'b/16$.
From (4), and because $u'B_iu=1-(u'x_i)^2$, we have $b'b=\sum\lbrack 1-(u'x_i)^2\rbrack^2=\lbrack\sum(u'x_i)^4\rbrack-4$.
After expansion: $c'c=\lbrack 1-(u_1^4+u_2^4+u_3^4)\rbrack/2$.
The orthogonality constraint reduces to $\sum(u'x_i)^3=24u_1u_2u_3=0$.
Either $u_1$ or $u_2$ or $u_3$ is null. We solve the optimization problem
in respect to the unkown vector $(u_1^2,u_2^2,u_3^2)'$, which describes the boundary of
the equilateral triangle having the three canonical base vectors as extreme points.
The norm of this latter unknown vector is minimized when its extremity lies on the mid of any of the side of the equilateral triangle:
$u_1^2=1/2,u_2^2=1/2,u_3^2=0$, or $u_1^2=1/2,u_2^2=0,u_3^2=1/2$, or $u_1^2=0,u_2^2=1/2,u_3^2=1/2$.
Thus the unit solutions vectors $u$ are in the directions of the edges and the maximal radius is $3/2$.
Note that minimizing $c'c$ let to retrieve the canonical base vectors as solutions, reproving part (a).

The results of theorem 1 are in agreement with those in \cite{Brandenberg2006,Devillers2003,Maehara2010}
and with those for equifacial tetrahedra \cite{Brandenberg2004,Theobald2003},
although the proof given here differs from the previously published ones.
\vskip 0.5cm

Returning to the case of the general tetrahedron,
it is pointed out that we are seeking for an unknown direction rather
than for an unit vector, so that the normalization condition can be weakened.
The leading term $u'u$ in (6) permits to write $F$ as an homogeneous quartic
to be minimized under the orthogonality constraint $\gamma'u=0$.
When this orthogonality constraint is satisfied by some vector $u$,
it is still satisfied by any vector colinear to $u$, and when the gradient
of the homogeneous quartic vanishes for some vector $u$,
it vanishes again for any vector colinear to $u$.
So, any non null vector $u$ got during the Newton iterations can be renormalized without
affecting the convergence of the iterative process.
Such technique was successfully used in a minimal variance computation
encountered in the geometric docking problem \cite{Petitjean2002,Petitjean2004}.
Some other practical implementation details follow.
Random starting unit vectors are generated from an isotropic distribution.
At the iteration $k$, the Lagrangians values are computed to minimize the
norm of the full gradient of $F$.
Then the Newton step $s_k=-\alpha_k H_k^{-1} g_k$ is computed, where $g_k$ and $H_k$
are respectively the gradient of $F$ at $u_k$ and its Hessian at $u_k$,
and $\alpha_k$ is a positive real, which is usually set to $1$ when the
quadratic convergence is achieved (see also the discussion at the end of the paper).
However the resulting vector $u_k+s_k$ does not satisfy anymore to the constraints.
The orthogonality condition can be restored via solving the minimization problem

\parindent 2cm
\vskip 0.5cm
               $Min_{\lbrace s\rbrace} (s-s_k)'(s-s_k)$, subject to
               \hskip 0.2cm $g_k's<\eta$, and to
               \hskip 0.2cm $u_{k+1}\gamma_{k+1}=0$ \hfill (7)
\vskip 0.5cm
\parindent 0.5cm

The parameter $\eta$ in the first constraint is set to ensure that, using a first order approximation,
there is a sufficient decrease of $F$. For that, we should impose a negative value to $g_k's$,
and we simply set $\eta=g_k's_k$.
But solving (7) needs a costly extra work due to the orthogonality constraint.
Thus, we expand the expression of the orthogonality constraint and we approximate it to the first order
of the unknown $s$. The modified constrained minimization problem (7) is solvable analytically,
and the Newton step $s$ is computed at low cost.
Then, as mentioned above, the normalization constraint is restored via setting $\Vert u_{k+1}\Vert=1$.

Our implementation of the solver was programmed in f77 with double precision (i.e. 64 bits) floating numbers.
The termination criterion was $\vert cos(u,c)\vert\le 10^{-10}$ and
$\vert (F_{k+1}-F_k)\vert/(F_{k+1}+F_k) \le 10^{-10}$.
Alternate termination criterions may consider the norm of the gradient,
the norm of the Newton step $s$, etc.

Thousands of random tetrahedra were generated from various probability laws,
including quasi-flat and quasi-linear sets.
The implementation of the Newton method appeared to be effective
for most random initial unit vectors.
In fact, for most tetrahedra, the convergence was observed for all initial vectors.
In the worst situations, reducing significantly the length of the Newton step at each iteration
before restoring the orthogonality constraint and operating with 100-250 initial values suffice
to overcome the problem of convergence failure and to observe a significant
number of times the convergence. Anyway, there may be several
locally minimal radius cylinders enclosing a tetrahedra, and the upper bound
of the number of local minima has been shown to be 9 \cite{Devillers2003}.
So, we indeed needed to operate with a sufficient number of random initial unit vectors.

A typical example is the regular tetrahedron
$x_1'=(1,1,1)$, $x_2'=(1,-1,-1)$, $x_3'=(1,-1,1)$, $x_4'=(-1,-1,1)$.
Two locally minimal radii were found: $\rho=1.41421356$ and $\rho=1.50000000$.
Among the 100 performed minimizations, the smallest radius was found 32 times
and the other one was found 68 times. The numbers of iterations ranged from 4 to 35,
with a mean value around 9.7.
The axis of the cylinder associated to the smallest radius 
is in agreement with the one given in theorem 1.
We observed that the other radius corresponds to the six unit vectors in the six respective directions
defined by the edges of the tetrahedron, for which the calculated value is indeed $\rho=3/2$.

\parindent 0cm
\vskip 1cm
\textbf{5. The circumscribed cylinder problem}
\vskip 0.5cm
\parindent 0.5cm

Here we set $n=5$. The convex hull of the $5$ points is either a tetrahedron with one of the
$5$ points lying in its interior, or a trigonal bipyramid.
Falling in this latter configuration is a necessary condition of existence of the circumscribed cylinder,
but it does not suffice.
According to (4), we must find the zeroes of the system (8) of 6 unknowns $c$ and $u$

\parindent 2cm
\vskip 0.5cm
               $u'B_iu=2c'x_i$ , $i=1,...,5$; \hskip 0.2cm $u'c=0$; \hskip 0.2cm $u'u=1$ \hfill (8)
\vskip 0.5cm
\parindent 0.5cm

The system (8) contains only 6 independant equations because the 5 first ones are dependant (sum to zero).
Again, the radius is computed from $c$ and $u$ and the center $c=\gamma/2$ is computed according to (5):

\parindent 2cm
\vskip 0.5cm
              $\rho^2=\bar{\Delta}=Tr(V)-u'Vu+c'c$

               $\gamma=T^{-1}X'b$, \hskip 0.2cm $b_i=u'B_iu$
\vskip 0.5cm
\parindent 0.5cm

Thus, reporting the expression of $c$ in some linear combination $t$ of the 5 first equations of (8)
lead to a system of only 3 equations of 3 unknowns, i.e. the components of $u$.
We define the vector $\mathbf{1}$ as a vector having its 5 components equal to 1.
$X$ being centered and of full rank, there is only one direction orthogonal to the $4$-dimensional subspace
defined by the the three independant columns of $X$ and the vector $\mathbf{1}$, this latter being itself orthogonal
to each of the 3 columns of $X$, i.e. $X'\mathbf{1}=0$.
We set the linear combination $t$ to follow this unique direction,
and we conventionnally normalize $t't=n$, and a plus sign is arbitrarily attributed to the component of $t$
having the largest absolute value. In other words, the computation of $t$ needs that
we perform a Gram-Schmidt orthogonalization of the columns of $X$ and of the vector $\mathbf{1}$.
Alternatively, an adequate rotation can be made so that $X$ is positioned in its principal components set of axis
(i.e. $T=X'X$ is diagonal), so that the columns of $X$ and the vector $\mathbf{1}$ are all orthogonal,
thus providing an easy way to compute $t$.

Having $X't=0$ means that $\sum\limits_{i=1}^{i=5} x_i=0$.
Then we define the matrix $M=\sum\limits_{i=1}^{i=5} t_iB_i$, and we consider the system (9)
issued from the linear combination $t$ of the 5 first equations of (8):

\parindent 2cm
\vskip 0.5cm
               $u'Mu=0$; \hskip 0.2cm $u'c=0$; \hskip 0.2cm $u'u=1$ \hfill (9)
\vskip 0.5cm
\parindent 0.5cm

This system, where $c$ is a quadratic function of $u$, has three equations of three unknowns.
It is valid if and only if $M$ is not null.

\vskip 0.5cm
\textbf{Theorem 2.}
\textit{$X$ being full dimensional, the matrix M is null if and only if two points are identical}.
\vskip 0.5cm

\textit{Proof.}
We assume $M=\sum t_iB_i=0$. We have $t'\mathbf{1}=0$, and so, $\sum t_i(I\cdot Tr(V_i) - V_i)=0$.
Because $Tr(M)=0$, we have also $\sum t_i Tr(Vi)=0$ and thus $\sum t_i x_ix_i'=0$.
We define the diagonal matrix $\Theta$ with its diagonal element $\Theta_{i,i}=t_i$, $i=1,...,5$.
Then, $X'\Theta X=0$, so that we can write the matrix product
\vskip 0.5cm

\parindent 2cm
$ \left(
  \begin{array}{c}
  X' \\
  \hline
  X'\Theta' \\
  \end{array}
  \right)
  \cdot
  \left(
  \begin{array}{c|c}
  X & \Theta X \\
  \end{array}
  \right)
  =
  \left(
  \begin{array}{c|c}
  T & 0 \\
  \hline
  0 & X'\Theta'\Theta X \\
  \end{array}
  \right)
$
\parindent 0.5cm

\vskip 0.5cm
The result of the matrix product above is a square matrix which cannot be of full rank because
it is the product of a 6 columns and 5 lines matrix by its transposed. Because $X$ is of full rank,
the block $X'\Theta'\Theta X=\sum t_i^2 x_ixi'$ must have at least one null eigenvalue.
So, there is some vector $\omega$ such that $\omega'X'\Theta'\Theta X\omega=0$,
and so $\Theta X\omega=0$.
We set $X_{1234}$ as the four lines array of the points $x_1$, $x_2$, $x_3$, $x_4$ ,
and $\Theta_{1234}$ as the associated diagonal submatrix of $\Theta$, and we remark that
$X_{1234}$ is of full rank because if it was not, there would be some vector orthogonal
to $x_1$, $x_2$, $x_3$ and $x_4$, and thus it would be
orthogonal to $x_5$ because $x_5=-\sum\limits_{i=1}^{i=4} x_i$,
thus implying that $X$ is not of full rank, a contradiction.
This remark stands for all the submatrices of $X$ containing four of the five points.
Now, if $\Theta_{1234}$ would be invertible, we would have found a vector $\omega$ such that $X_{1234}\omega=0$,
which is impossible, so $\Theta_{1234}$ is not invertible, and there at least one zero element
among $t_1$, $t_2$, $t_3$, $t_4$.
This conclusion is valid for all four-tuples of $\lbrace 1,2,3,4,5\rbrace$,
so that there are at least two distinct null elements among $t_1$, $t_2$, $t_3$, $t_4, t_5$.

Conventionnally, we set $t_4=0$ and $t_5=0$ without loss of generality.
From the orthogonality conditions $X't=0$ and $\textbf{1}'t=0$, we have
$t_1x_1+t_2x_2+t_3x_3=0$ and $t_1+t_2+t_3=0$, which means that $x_1$, $x_2$, $x_3$ are aligned.
Assuming $t_1t_2t_3\ne 0$, $x_3=(t_1x_1+t_2x_2)/(t_1+t_2)$ and expanding the expression of
$M=t_1x_1x_1'+t_2x_2x_2'+t_3x_3x_3'$ leads to $M=t_1t_2(x_1-x_2)(x_1-x_2)'/(t_1+t_2)$.
The assumption $M=0$ would imply here that $x_1=x_2=x_3$,
which is impossible due to the full dimensionnality of $X$.
So, $t_1t_2t_3=0$, and we must assume that a third element of $t$ is null.
We set conventionnally $t_3=0$, and we have both $t_1+t_2=0$ and $t_1x_1+t_2x_2=0$,
meaning that $x_1=x_2$ because $t$ is not null by definition, and $M$ is indeed null.

Conversely, we consider a full rank set such that two points,
say $1$ and $2$, are identical.
The unique direction $t$ is colinear to $(1,-1,0,0,0)$ and thus $M=0$.

\vskip 0.5cm
Having two identical points means that in fact we face to the four points problem,
and this situation was considered in section 4. It may also be considered that suppressing
the equation $u'Mu=0$ lead to an underdetermined system, and so the problem of finding
the zeroes of a function should be transformed in an optimization problem
where the radius of the cylinder is to be minimized.
Anyway, we further assume that we have never two identical points.

\vskip 0.5cm
\textbf{Theorem 3.}
\textit{The set of solutions of (8) and the set of solutions of (9) are equal}.
\vskip 0.5cm

\textit{Proof.}
Obviously, the set of solutions of (8) is included in the set of solutions of (9).
In order to ensure that the solutions of (9) satisfies to (8), we should check
that they indeed satisfy to the five equations $u'B_iu=2c'x_i$, $i=1,...,5$,
or, in matricial form, $b=2Xc$.
We recall that the center is computed from $\gamma=2c=T^{-1}X'b$, with
$b_i=u'B_iu$, $i=1,...,5$. The $5$-dimensional vector $b$ can be decomposed
in the basis of $5$ independant vectors build with the $3$ columns of $X$
and the vector $\textbf{1}$ and the vector $t$, i.e.
$b=2X\beta_X +\beta_1 \textbf{1} +\beta_t t$, where $\beta_X$ contains
the three coefficients associated to the columns of $X$ and $\beta_1$ and $\beta_t$
are the coefficients associated respectively to $\textbf{1}$ and $t$.
The sum of the five $B_i$ is null, thus $1'b=0$ and so $\beta_1=0$.
From $u'Mu=0$ we have $t'b=\sum t_i(u'B_iu)=0$.
Thus, $2t'X\beta_X +\beta_t t't=0$, from which $\beta_t=0$ because $t'X=0$.
So, $b$ is indeed a linear combination of the columns of $X$, i.e. $b=2X\beta_X$.
The computed center is $2c=T^{-1}X'b$, so $\beta_X=c$ and indeed $b=2Xc$.
Remark: the assumption $M\ne 0$ was not used in the proof above.

\vskip 0.5cm
\textbf{Lemma 1.}
\textit{There are at most 6 cylinders circumscribed to 5 points}.
\vskip 0.5cm

\textit{Proof.}
The homogeneous polynomial equations $u'Mu=0$ and $u'c=0$ in the system (9)
define the intersection of two curves in the projective plane of respective degrees 2 and 3.
The normalization condition $u'u=1$ is just a practical way to avoid a null vector $u$,
and in fact all non null vectors colinear to a solution $u$ at the intersection of $u'Mu=0$ and $u'c=0$
define the same intersection and lead to the same cylinder.
So, applying Bezout theorem shows that we should have at most 6 cylinders passing through 5 points.
This result was already proved in \cite{Devillers2003} via an other method.

\vskip 0.5cm
\textbf{Lemma 2.}
\textit{There is no cylinder circumscribed to 5 points when
either $M$ or $-M$ is positive definite}.
\vskip 0.5cm

\textit{Proof.}
Obvious from $u'Mu=0$.

\vskip 0.5cm
The matrix $M=\sum t_i(I\cdot Tr(V_i)-V_i)=I\cdot Tr(X'\Theta X)-X'\Theta X$,
characterizes a quadric surface.
Except in the particular situation where an eigenvalue of $M$ vanishes,
the surface $u'Mu=t'b=0$, i.e. $\sum t_i(u'u\cdot x_i'x_i-(u'x_i)^2)=0$,
is an elliptical cone.
The cubic surface is $u'\gamma=0$,
with $\gamma=T^{-1}\sum x_i(u'u\cdot x_i'x_i-(u'x_i)^2)$.

\vskip 0.5cm
Solving the system (9) can be reduced to an one unknown problem.
We assume without loss of generality
that the data set is rotated to have $M$ being diagonal with eigenvalues
$\mu_1\ge\mu_2\ge\mu_3$.
The case $\mu_1\mu_3>0$ was proved to correspond
to no circumscribed cylinder and the case $\mu_1=\mu_3=0$ was proved to correspond
to an infinite number of cylinders (four points problem).

\vskip 0.5cm
The equations $u'Mu=0$ and $u'u=1$ define a linear system of two equations
of the the three squared components of $u$.
One of the components $u_1$ or $u_2$ or $u_3$ is taken as the unknown parameter,
and the expressions of the two other components are reported in $u'\gamma=0$.
Due to the free choice of the signs of these two other components,
we have in fact four one unknown equations to solve.
Assuming that we performed the rotation above and that the eigenvalues
of $M$ are separated, the relations between the components of $u$ are

\parindent 2cm
\vskip 0.5cm
  $u_2^2=(u_1^2(\mu_3-\mu_1)-\mu_3)/(\mu_2-\mu_3) \hskip 0.5cm
   u_3^2=(u_1^2(\mu_1-\mu_2)+\mu_2)/(\mu_2-\mu_3)$

  $u_1^2=(u_2^2(\mu_3-\mu_2)-\mu_3)/(\mu_1-\mu_3) \hskip 0.5cm
   u_3^2=(u_2^2(\mu_2-\mu_1)+\mu_1)/(\mu_1-\mu_3)$

  $u_1^2=(u_3^2(\mu_2-\mu_3)-\mu_2)/(\mu_1-\mu_2) \hskip 0.5cm
   u_2^2=(u_3^2(\mu_3-\mu_1)+\mu_1)/(\mu_1-\mu_2)$
\vskip 0.5cm
\parindent 0.5cm

The components of $u$ must be in $\lbrack 0;1\rbrack$, but from the relations above
we get tighter bounds for the unknown component to be selected. When $\mu_2\ge 0$:

\parindent 2cm
\vskip 0.5cm
                      $0 \le u_1^2\le 1-\mu_1/(\mu_1-\mu_3)$

                      $0 \le u_2^2\le 1-\mu_2/(\mu_2-\mu_3)$

   $ \mu_2/(\mu_2-\mu_3) \le u_3^2\le 1+\mu_3/(\mu_1-\mu_3)$
\vskip 0.5cm
\parindent 0.5cm

And when $\mu_2\le 0$:

\parindent 2cm
\vskip 0.5cm
   $ -\mu_2/(\mu_1-\mu_2) \le u_1^2\le 1-\mu_1/(\mu_1-\mu_3)$

                       $0 \le u_2^2\le 1+\mu_2/(\mu_1-\mu_2)$

   $                    0 \le u_3^2\le 1+\mu_3/(\mu_1-\mu_3)$
\vskip 0.5cm
\parindent 0.5cm

When either $\mu_1=\mu_2$ or $\mu_2=\mu_3$, some of the relations
and intervals above are invalid.
The case $\mu_1=\mu_3$ is such than there is no cylinder to find,
either because $M$ or $-M$ is positive definite, or because $M=0$
and we have a four points problem.
Outside these situations, we can always retain $u_2$ as the unknown parameter.
Whatever unknown we select, either $u_1$, $u_2$ or $u_3$,
the resulting function is not a polynomial.

\vskip 0.5cm
\textbf{Theorem 4.}
\textit{The system (9) can been solved via extracting the real roots
of a polynomial of degree at most equal to 6.}
\vskip 0.5cm

\textit{Proof.}
Having rotated the set so that $M$ is diagonal, we rotate it again
around the second eigenvector of $M$ with an angle $\alpha_2$ to be specified.

\parindent 1cm
\vskip 0.5cm
 $u'Mu = u_1^2(\mu_1cos^2(\alpha_2)+\mu_3sin^2(\alpha_2))
       + u_2^2\mu_2
       + u_3^2(\mu_1sin^2(\alpha_2)+\mu_3cos^2(\alpha_2))$

\parindent 2.5cm
       $- 2u_1u_3sin(\alpha_2)cos(\alpha_2)(\mu_1-\mu_3)$
\vskip 0.5cm
\parindent 0.5cm

We consider the general case where $\mu_1\mu_3<0$.
The case where $\mu_1\mu_3=0$ and either $\mu_1\ne 0$ or $\mu_3\ne 0$
will be considered later.

We select $\alpha_2$ so that $\mu_1sin^2(\alpha_2)+\mu_3cos^2(\alpha_2)=0$,
i.e. $tg(w)=\sqrt{-\mu_3/\mu_1}$.
The quadric becomes a parabola in the projective plane
defined by the coordinate $u_1$.

\parindent 2cm
\vskip 0.5cm
     $u_1^2(\mu_1+\mu_3)+u_2^2\mu_2-2u_1u_3\sqrt{-\mu_1\mu_3}=0$ \hfill (10)

     $(\frac{u_3}{u_1}) =
      - \frac{\mu_2}{2\sqrt{-\mu_1\mu_3}} (\frac{u_2}{u_1})^2
      - \frac{\mu_1+\mu_3}{2\sqrt{-\mu_1\mu_3}}$ \hfill (11)

\vskip 0.5cm
\parindent 0.5cm

Reporting $u_3/u_1$ in the expression of the cubic curve in the projective plane
leads to a polynomial of degree six in $u_2/u_1$, which has at most 6 real roots.
We deduce $u_3/u_1$ from (11) and then we have the direction of $u$.
So we get at most 6 cylinders.

The existence of solutions such that $u_1=0$ should be checked.
When $\mu_2\ne 0$, we deduce from (10) that $u_2=0$, and hence $u'=(0,0,\pm 1)$
should satisfy to $u'\gamma=0$, so that the coefficient of $u_3^3$ in the cubic form
is null. If the solution $u'=(0,0,\pm 1)$ is indeed found, the other ones
found in the projective plane from (11) are such that the curve defined from
$u'\gamma=0$ degenerates to a quintic of $u_2/u_1$
and we cannot expect more than 5 real roots, i.e. there are at most 6 cylinders.

When $\mu_2=0$, the homogeneous cubic of $u_2$ and $u_3$ offers
at most 3 real roots when solved either in $u_2/u_3$ or $u_3/u_2$
($u_2$ and $u_3$ cannot be both null).
The solutions other than $u'=(0,u_2,u_3)$
found in the projective plane from (11) are such that the curve defined from
$u'\gamma=0$ is a cubic of $u_2/u_1$ and we expect at most 3 real roots,
i.e. there are at most 6 cylinders.

We consider now the case where $\mu_1\mu_3=0$ and either $\mu_1\ne 0$ or $\mu_3\ne 0$.
Selecting $\alpha_2=0$ so that $M$ is diagonal, we have
$\mu_1u_1^2+\mu_2u_2^2+\mu_3u_3^2=0$.
When $\mu_2\ne 0$, there is only one potential solution which is either
$u'=(0,0,\pm 1)$ or $u'=(\pm 1,0,0)$.
When $\mu_2=0$, either $u_1=0$ when $\mu_1\ne 0$ or $u_3=0$ when $\mu_3\ne 0$.
When $\mu_1\ne 0$, $u'\gamma$ is an homogeneous cubic of $u_2$ and $u_3$
which offers at most 3 roots when solved either in $u_2/u_3$ or $u_3/u_2$
($u_2$ and $u_3$ cannot be both null), so that we expect at most 3 cylinders.
A similar conclusion is got when $\mu_3\ne 0$.

In any case, the system (9) can be solved
by finding the roots of a polynomial of degree at most 6,
and we cannot get more than 6 cylinders, which redemonstrates lemma 1.

\parindent 0cm
\vskip 1cm
\textbf{6. An example of trigonal bipyramid with a symmetry axis of order 3}
\vskip 0.5cm
\parindent 0.5cm

We present the case of a bipyramid $\cal{BP}$ symmetric around its equilateral triangular basis
with a symmetry axis of order 3 orthogonal to the triangular basis.
We set
$x_1'=(0,0,h)$, $x_2'=(0,0,-h)$, $x_3'=(1,0,0)$, $x_4'=(-1/2,\sqrt{3}/2,0)$, $x_5'=(-1/2,-\sqrt{3}/2,0)$.
The set is centered.
\parindent 0cm

\vskip 0.5cm
   $
     I\cdot Tr(V_1)-V_1 =
     I\cdot Tr(V_2)-V_2 = \left(
     \begin{array}{ccc}
     h^2 & 0   & 0 \\
     0   & h^2 & 0 \\
     0   & 0   & 0 \\
     \end{array}
         \right)
   \hskip 0.5cm
     I\cdot Tr(V_3)-V_3 = \left(
     \begin{array}{ccc}
     0 & 0 & 0 \\
     0 & 1 & 0 \\
     0 & 0 & 1 \\
     \end{array}
         \right)
   $

\vskip 0.25cm
   $
     I\cdot Tr(V_4)-V_4 = \left(
     \begin{array}{ccc}
     3/4         &  \sqrt{3}/4 & 0 \\
      \sqrt{3}/4 & 1/4         & 0 \\
     0           & 0           & 1 \\
     \end{array}
         \right)
   \hskip 0.5cm
     I\cdot Tr(V_5)-V_5 = \left(
     \begin{array}{ccc}
     3/4         & -\sqrt{3}/4 & 0 \\
     -\sqrt{3}/4 & 1/4         & 0 \\
     0           & 0           & 1 \\
     \end{array}
         \right)
   $

\vskip 0.25cm
   $
     t = (1/\sqrt{6}) \left(
     \begin{array}{c}
      3 \\
      3 \\
     -2 \\
     -2 \\
     -2 \\
     \end{array}
         \right)
   \hskip 1.5cm
     M = (1/\sqrt{6}) \left(
     \begin{array}{ccc}
     6h^2-3 &        &  0 \\
     0      & 6h^2-3 &  0 \\
     0      & 0      & -6 \\
     \end{array}
         \right)
   $

\vskip 0.5cm
So, there is no circumscribed cylinder when $h< \sqrt{2}/2$.
We assume $h\ge \sqrt{2}/2$.

\vskip 0.5cm
   $
     B_1 = B_2 = \left(
     \begin{array}{ccc}
     (6h^2-3)/10 &           0 & 0    \\
     0           & (6h^2-3)/10 & 0    \\
     0           &           0 & -3/5 \\
     \end{array}
         \right)
   $

\vskip 0.25cm
   $
     B_3 = \left(
     \begin{array}{ccc}
     (-4h^2-3)/10 & 0           & 0   \\
     0            & (7-4h^2)/10 & 0   \\
     0            & 0           & 2/5 \\
     \end{array}
         \right)
   $

\vskip 0.25cm
   $
     B_4 = \left(
     \begin{array}{ccc}
     (9-8h^2)/20 &   \sqrt{3}/4 & 0   \\
      \sqrt{3}/4 & (-1-8h^2)/20 & 0   \\
               0 &            0 & 2/5 \\
     \end{array}
         \right)
   $

\vskip 0.25cm
   $
     B_5 = \left(
     \begin{array}{ccc}
     (9-8h^2)/20 &  -\sqrt{3}/4 & 0   \\
     -\sqrt{3}/4 & (-1-8h^2)/20 & 0   \\
               0 &            0 & 2/5 \\
     \end{array}
         \right)
   $
\vskip 0.5cm

   $b_1=b_2=3(2h^2-2+u_1^2+u_2^2)/10$

   $b_3=(4-u_1^2(4h^2+7)+u_2^2(3-4h^2))/10$

   $b_4=(8+(1-8h^2)u_1^2-(9+8h^2)u_2^2+10\sqrt{3}u_1u_2))/20$

   $b_5=(8+(1-8h^2)u_1^2-(9+8h^2)u_2^2-10\sqrt{3}u_1u_2))/20$


\vskip 0.25cm
   $
     T^{-1}X' = \left(
     \begin{array}{ccccc}
     0          & 0           & 2/3 & -1/3       & -1/3        \\
     0          & 0           & 0   & \sqrt{3}/3 & -\sqrt{3}/3 \\
     1/2h       & -1/2h       & 0   & 0          & 0           \\
     \end{array}
         \right)
   \hskip 0.5cm
     \gamma=2c = \left(
     \begin{array}{c}
     (-u_1^2+u_2^2)/2 \\
     u_1u_2           \\
     0                \\
     \end{array}
         \right)
   $

\vskip 0.5cm
Solving $u'Mu$=0 with $u'\gamma=0$, i.e.
$(2h^2+1)(u_1^2+u_2^2)=2$ and $u_1(u_1^2-3u_2^2)=0$
leads to the six desired directions $u$ and we easily deduce Theorem 5 below.

\vskip 0.5cm
\textbf{Theorem 5.}
\textit{\\(a) The symmetric bipyramid $\cal{BP}$ of interapex half-height $h$ has no circumscribed cylinder when $h<\sqrt{2}/2$.
\\(b) When $h\ge\sqrt{2}/2$, $\cal{BP}$ has six circumscribed cylinders.
The axis are in the directions $u'=(0,\sqrt{2/(2h^2+1)},\pm\sqrt{(2h^2-1)/(2h^2+1}))$,
which intersect at $c'=(1/(4h^2+2),0,0)$,
plus their rotated images of respective angles $2\pi/3$ and $4\pi/3$ around the interapex symmetry axis of order 3.
All six cylinders have radius $\rho=(4h^2+1)/(4h^2+2)$.
}

\vskip 0.5cm
Neither the condition of existence of the cylinders circumscribed to $\cal{BP}$
nor their analytical calculations seem to be previously reported in the literature.
They are mainly consequences of Theorem 3 and Lemma 2, which themselves cannot
be trivially deduced from existing results in cited papers and provide
a simpler approach to the circumscribed cylinder problem in $E^3$.

\parindent 0cm
\vskip 1cm
\textbf{7. Discussion and conclusion}
\vskip 0.5cm
\parindent 0.5cm

Some results in this work were already available in previous papers such as
\cite{Brandenberg2004,Devillers2003,Theobald2003}.
We had clearly mentioned these references at the appropriate places.
Since these results are obtained here via different techniques and add both
to the clarity and the self-content character of the paper, we left them.
Having said that, we recall that several novelties are introduced in this paper.
First, formulating the problems through variance minimizations led us
in the three-dimensional case to the original systems (6) and (9),
which offer particularly simple expressions of the unknown axis $u$.
It means that iterative numeric solvers need only random $3$-tuples,
and can be programmed by most users having some knowledge
of constrained optimization, without needing to install specialized
mathematical packages.

In the case of the circumscribed cylinder to five points,
we produced a simple and new condition of existence of the cylinder
relying on the eigenvalues of the matrix $M$ in section 5.
This matrix and its properties were yet unpublished.
The explicit expressions of the polynomials of one variable
at the end of section 5 are original, too.
However, we found easier to programme directly the solver of (9)
rather than seeking for the roots of the one variable polynomials.
It was just an arbitrary developper's choice.
Then, the third reviewer of this paper mentioned that
the presence of square roots in the polynomials could be penalizing,
and thus particular methods such as \cite{Eigen2005} could be used.

Finally, the original analytical solution for the symmetric
bipyramid in section 6, which relies on the explicit calculation of
the matrix $M$, is outlined.

By no way we claim that the original present approach
is better than the previous ones, such as in \cite{Devillers2003},
and the present software, named CYL, downloadable for free at
\href{http://petitjeanmichel.free.fr/itoweb.petitjean.freeware.html}{http://petitjeanmichel.free.fr/itoweb.petitjean.freeware.html},
is neither claimed to offer the best possible implementations of the methods
nor is claimed to work better than other cylinder computation softwares.
Moreover, the second reviewer informed us that the algebraic approach in \cite{Devillers2003}
could guarantee the result, although our solver (see in particular section 4) cannot offer this guarantee
for all datasets, despite the encouraging results we got.

It would have been of high interest to compare our practical experiment results
with those of existing methods such as the algebraic method in \cite{Devillers2003}.
But, apart for the regular tetrahedron (detailed numerical results
at the end of section 5), and for the bipyramid in section 6,
we found only one public dataset, mentioned in \cite{Watson2006}:
we retrieved the optimal radius of the 12 points set minimal enclosing cylinder
computed in \cite{Watson2006} with all significant digits.

We propose to those readers who have access to some cylinders computation
software, to exchange with us any data and practical experimentation
results they like.

\parindent 0cm
\vskip 1cm
\textbf{Acknowledgements}
\vskip 0.5cm
\parindent 0.5cm
I thank the three reviewers for having taken time to analyse this work,
and I am particularly grateful to the third reviewer who did pertinent
remarks about computational details and provided ref. \cite{Eigen2005}.

\vskip 0.5cm
\parindent 0cm
\renewcommand{\refname}{\small{References}}

\end{document}